\documentstyle{article}
\begin{document}
\begin{center}
{ \bf \Large Improved Results on Robust Stability of Multivariable
Interval Control Systems \footnote[1] {\small Supported by
National Key Project and National Natural Science Foundation of
China(69925307). } }
\end{center}

\vskip 0.6cm


\centerline{\large Zhizhen Wang$^{(1)}$ \ \ \ Long Wang$^{(2)}$ \
\ \ Wensheng Yu$^{(1)}$ } \vskip 6pt \centerline{(1)\ \ Institute
of Automation, Chinese Academy of Sciences, } \centerline{Beijing,
100080, CHINA} \centerline{(2)\ \ Center for Systems and Control,
Department of Mechanics and Engineering Science, }
\centerline{Peking University, Beijing 100871, CHINA }

\vskip 0.6cm

\vskip 6pt
\begin{minipage}[t]{14cm}{
{\bf Abstract: }
For  interval polynomial matrices, we identify the minimal testing set, whose
stability can guarantee that of the whole uncertain set. Our results improve the conclusions given by Kamal and Dahleh.

\vskip 4pt

 }
\end{minipage}

\section{Introduction}
\par \indent

A control system is said to be robust if it retains an assigned
degree of stability or performance under perturbations. Robustness
is considered as an elementary aspect in the analysis and design of
control systems. There are several research lines for different
models of uncertainties. Robust control under structured perturbations is an
active area in recent years$^{\cite{khar}-\cite{barm1}}$.
\par\indent

The parametric robust stability problem was initiated by
Kharitonov(1978)$^{\cite{khar}}$, who gave a sufficient and necessary
condition for the robust stability of interval polynomials. Chapellat et al.(1989)$^{\cite{chap1}}$ proved the Box Theorem
for the multilinear model $P(s)=U(s)V(s)+X(s)Y(s)$.
\par\indent

When considering multivariable systems, the
robust stability problem becomes more complicated$^{\cite{kama,long,chap1,bias,barm,koka}}$.
Ref. \cite{bias} and Ref.
\cite{barm} discussed the scalar matrix case and pointed out that a
finite test does not hold. Ref. \cite{long}, \cite{koka}  presented an
$n$-dimensional result for polynomial matrices.  Ref. \cite{kama} gave a
$2n$-dimensional result for the model $B(s)A(s)+D(s)C(s)$, where
$B(s),D(s)$ are interval polynomial matrices and $A(s),C(s)$ are
fixed ones.
\par\indent

The purpose of this paper is to improve the robust stability criteria for multivariable interval
control systems obtained by Kamal and Dehleh. By exploiting the uncertainty structures, we are able to reduce
the computational burden in checking robust stability.
\section{Preliminaries}
{\bf Definition 2.1}\ \ Given an interval polynomial set
$r(s)=\sum_{i=0}^m r_i s^i, \ \  r_i \in [r_i^L, r_i^U] $, its
Kharitonov vertex set is
 $r_V(s)=\{r_1(s), r_2(s), r_3(s), r_4(s)\}$, and its Kharitonov edge set is
 $r_E(s)=\{\lambda r_i(s)+(1-\lambda) r_j(s), \ \
(i,j)\in \{(1,2), (2,4), (4,3), (3,1)\} \ \  \lambda \in [0,1]\}$,
 where
  $$
  \begin{array}{l}
  r_1(s)=r_0^L+r_1^L s+r_2^U s^2+r_3^U s^3+r_4^L s^4+\dots\\
  r_2(s)=r_0^L+r_1^U s+r_2^U s^2+r_3^L s^3+r_4^L s^4+\dots\\
  r_3(s)=r_0^U+r_1^L s+r_2^L s^2+r_3^U s^3+r_4^U s^4+\dots\\
  r_4(s)=r_0^U+r_1^U s+r_2^L s^2+r_3^L s^3+r_4^U s^4+\dots
  \end{array}
  $$

{\bf Definition 2.2}\ \ A polynomial is stable if all its zeroes belong to
the open left half plane. A polynomial matrix is stable if its determinant is stable.

{\bf Lemma 1}\quad (Box Theorem$^{\cite{chap1}}$) Consider the
polynomial family $\Delta(s)=\{\delta(s,p)=F_1(s)
P_1(s)+\dots+F_m(s) P_m(s), \; P_i(s)\mbox{ are interval
polynomials}, F_i(s) \mbox{ are fixed polynomials, } i=1, \dots, m
\}$ and suppose $\Delta(s)$ is without degree dropping. Then
$\Delta(s)$ is Hurwitz stable if and only if $\Delta_E (s)$ is
Hurwitz stable, where $\Delta_E (s)=\cup_{l=1}^m
\{\sum_{i=1}^{l-1}F_i (s)K_{P_i}(s)+ F_l(s)
E_{P_l}(s)+\sum_{i=l+1}^m F_i (s)K_{P_i}(s)\}$ and
$K_{P_i}(s),E_{P_i}(s)$ are the Kharitonov vertex set, Kharitonov
edge set of $P_i(s)$, respectively. (Here we take $\sum_{i=r}^t
f_i=0$ if $r>t$).

Let $A(s)=\left(a_{ij}(s)\right)_{n\times n},C(s)=\left(c_{ij}(s)\right)_{n\times n}$
be two fixed matrices, where $A(s)=\left(a_{ij}(s)\right)_{n\times n}$ means that the matrix is an $n\times n$ one
and its $ij$-th entry is $a_{ij}(s)$.   Given two interval polynomial matrices
$B(s)=\left(b_{ij}(s)\right)_{n\times n},D(s)=\left(d_{ij}(s)\right)_{n\times n}$, i.e.,
$b_{ij}(s),d_{ij}(s)$ are interval polynomials. Let
\begin{eqnarray}
B(ij)=\left\{b_{ij}(s):b_{ij}(s)=b_{ij}^{(0)}+b_{ij}^{(1)}s+\dots+b_{ij}^{(n)} s^n,\;
b_{ij}^{(k)}\in[\underline b_{ij}^{(k)},\overline b_{ij}^{(k)}]\right\}\\
D(ij)=\left\{d_{ij}(s):d_{ij}(s)=d_{ij}^{(0)}+d_{ij}^{(1)}s+\dots+d_{ij}^{(n)} s^n,\;
d_{ij}^{(k)}\in[\underline d_{ij}^{(k)},\overline d_{ij}^{(k)}]\right\}\\
{\cal B}(s)=\left\{B(s):B(s)=\left(b_{ij}(s)\right)_{n\times n},\; b_{ij}(s)\in B(ij)\right\}\\
{\cal D}(s)=\left\{D(s):D(s)=\left(d_{ij}(s)\right)_{n\times n},\; d_{ij}(s)\in D(ij)\right\}
\end{eqnarray}
Denote  Kharitonov vertex sets and Kharitonov edge sets of
$B(ij),D(ij)$ as ${\cal B}_V(ij),{\cal D}_V(ij),{\cal B}_E(ij)$
and ${\cal D}_E(ij)$ respectively.

\section{ Main Results}
\par \indent

In this part, we consider the uncertain family
\begin{equation}\label{1}
{\cal M}(s)={\cal B}(s)A(s)+{\cal D}(s)C(s)
\end{equation}

Let $S_n$ denote  the set of all bijections of the set
$\{1,\dots,n\}$ onto itself. The following theorem was given by
Kamal et al.(1996).

{\bf Proposition 1}$^{\cite{koka}}$\ \ ${\cal M}(s)$ is stable for
all $B(s)\in{\cal B}(s),D(s)\in{\cal D}(s)$ if and only if ${\cal
M}(s)$ is stable for all $B(s)\in {\cal B}_E(s),D(s)\in{\cal
D}_E(s)$, where
\begin{eqnarray}
{\cal B}_E(s)=\left\{B(s)\in{\cal B}(s):\; \begin{array}{l}
b_{ij}(s)\in {\cal B}_E(ij), \; j=\sigma(i)\\
b_{ij}(s)\in {\cal B}_V(ij), \; j\not=\sigma(i)\\
\sigma\in S_n
\end{array}\right\}\\
{\cal D}_E(s)=\left\{D(s)\in{\cal D}(s):\; \begin{array}{l}
d_{ij}(s)\in {\cal D}_E(ij), \; j=\sigma(i)\\
d_{ij}(s)\in {\cal D}_V(ij), \; j\not=\sigma(i)\\
\sigma\in S_n
\end{array}\right\}
\end{eqnarray}

As a generalization of  the Generalized Kharitonov Theorem, Proposition 1
addresses the robust stability of ${\cal M}(s)$ and reduce it
to $4^{2n^2}(n!)^2$ robust stability problems involving
determinants with $2n$ parameters. In this paper, we
will show that, by making use of the uncertainty structure information, the
original robust stability problem can be reduced to $\sum_{t=0}^n (C_n^t)^2 4^{2n^2}(n!)$
robust stability problems involving $n$ parameters, where $C_n^t=\frac{n!}{t!(n-t)!}$. In the
sequel, the notation $X\times Y$ stands for the Cartesian product, which
is the set of all ordered pairs $(x,y)$, where $x\in X,y\in Y$.
The following lemma is needed before we present our main result.

{\bf Lemma 2}\ \  $A(s), B(s),C(s)$ and $D(s)$  are defined as before. Let
\begin{eqnarray}
N^{(i)}(s)=\left\{B(s)A(s)+D(s)C(s):\; \begin{array}{l}
\exists k_1,k_2\in\{1,\dots,n\}; l=1,\dots,n\\
b_{ik_1}\in {\cal B}_E(ik_1); b_{il}\in {\cal B}_V(il),\; l\not=k_1\\
d_{ik_1}\in {\cal D}_E(ik_2); d_{il}\in {\cal D}_V(il),\;
l\not=k_2
\end{array}\right\}\\
N^{(i)}_E(s)=\left\{B(s)A(s)+D(s)C(s):\; \begin{array}{l}
\exists k\in\{1,\dots,n\}; l=1,\dots,n\\
b_{ik}\times d_{ik}\in  ({\cal B}_E(ik)\times {\cal D}_V(ik))\cup ({\cal B}_V(ik)\times {\cal D}_E(ik))\\
b_{il}\times d_{il}\in  ({\cal B}_V(il)\times {\cal D}_V(il)),
l\not=k
\end{array}\right\}
\end{eqnarray}
Then, for all $i\in\{1,\dots,n\},
N^{(i)}(s) \mbox{ is stable if and only if } N^{(i)}_E(s)\mbox{ is stable }.$\\
Proof\ \ (Sufficiency) Suppose that $N^{(i)}_E(s)\mbox{ is stable }$, our aim is to prove that
$N^{(i)}(s) \mbox{ is stable}$.
For any $N(s)\in N^{(i)}(s)$, there exist $k_1,k_2\in\{1,\dots,n\}$ such that
\begin{equation}\label{10}
 \begin{array}{ll}
 b_{ik_1}\in {\cal B}_E(ik_1);& b_{il}\in {\cal B}_V(il),\; l\not=k_1,l=1,\dots,n\\
d_{ik_1}\in {\cal D}_E(ik_2);& d_{il}\in {\cal D}_V(il),\;
l\not=k_2,l=1,\dots,n
\end{array}
\end{equation}
and
\begin{equation}\label{11}
N(s)=\left(
\begin{array}{lll}
\sum_k \left(b_{1k}(s)a_{k1}(s)+d_{1k}(s)c_{k1}(s)\right)&\dots&\sum_k \left(b_{1k}(s)a_{kn}(s)+d_{1k}(s)c_{kn}(s)\right)\\
\dots&\dots&\dots\\
\sum_k \left(b_{ik}(s)a_{k1}(s)+d_{ik}(s)c_{k1}(s)\right)&\dots&\sum_k \left(b_{ik}(s)a_{kn}(s)+d_{ik}(s)c_{kn}(s)\right)\\
\dots&\dots&\dots\\
\sum_k \left(b_{nk}(s)a_{k1}(s)+d_{nk}(s)c_{k1}(s)\right)&\dots&\sum_k \left(b_{nk}(s)a_{kn}(s)+d_{nk}(s)c_{kn}(s)\right)
\end{array}\right)
\end{equation}
By definition, $N(s)$ is stable if and only if ${\mathrm{det}}N(s)$ is stable. By Laplace formula, expanding the determinant
along the $i$-th row, we have
$$
{\mathrm{det}}N(s)=\sum_k \left(b_{ik}(s)a_{k1}(s)+d_{ik}(s)c_{k1}(s)\right)N_{i1}+\dots+
\sum_k \left(b_{ik}(s)a_{kn}(s)+d_{ik}(s)c_{kn}(s)\right)N_{in},
$$
where $N_{ij}$ is the algebraic complementary minor of the $ij$-th entry of $N(s)$, which is independent of $b_{i1}(s),\dots,
b_{in}(s),d_{i1}(s),\dots,d_{in}(s)$. By a simple manipulation, we have
$$
\begin{array}{lll}
{\mathrm{det}}N(s)&=&\sum_{j=1}^n b_{ij}(s)\left(\sum_{k=1}^n
a_{jk}(s)N_{ik}\right)+
\sum_{j=1}^n d_{ij}(s)\left(\sum_{k=1}^n c_{jk}(s)N_{ik}\right)\\
&=& b_{ik_1}(s)\delta_1(s)+d_{ik_2}(s)\delta_2(s)+\delta_3(s)
\end{array}
$$
where $\delta_l(s)$ is a term which is independent
of $b_{i1}(s),\dots, b_{in}(s),d_{i1}(s),\dots,d_{in}(s)$. Denote
\begin{equation}\label{4}
\left\{\begin{array}{l} b_{ik_1}(s)\times d_{ik_2}(s)\in
\left({\cal B}_E(i k_1)\times {\cal D}_V(i k_2)\right)\cup
\left({\cal B}_V(i k_1)\times {\cal D}_E(i k_2)\right)\\
b_{il}(s)\in {\cal B}_V(il),l\not=k_1,l=1,\dots,n\\
d_{il}(s)\in {\cal D}_V(il),l\not=k_2,l=1,\dots,n
\end{array}\right.
\end{equation}
By Lemma 1, ${\mathrm{det}}N(s)$ is stable for $N(s)$ satisfying
(\ref{10}),(\ref{11}) if and only if ${\mathrm{det}}N(s)$ is stable
for $N(s)$ satisfying (\ref{11}),(\ref{4}). Since $N^{(i)}_E(s)\mbox{ is stable }$,
it shows that
$N(s)$ is stable. Hence, sufficiency is proved. \\
(Necessary) This part is obvious, since $N^{(i)}_E(s)\subset
N^{(i)}(s)$.
\par\indent

Before presenting our main result, we first define some notations.
Let $T_n$ be the set of all functions from the set
$\{1,\dots,n\}$ into itself. Take
$$
{\cal BD}=\left\{B(s)\times D(s): \begin{array}{ll}
b_{ij}(s)\times d_{ij}(s)\in \left({\cal B}_E(ij)\times {\cal
D}_V(ij)\right)\cup \left({\cal B}_V(ij)\times
{\cal D}_E(ij)\right)& j=\eta(i)\\
b_{ij}(s)\times d_{ij}(s)\in {\cal B}_V(ij)\times {\cal D}_V(ij)
&j\not= \eta(i)\\
\eta\in T_n&
\end{array}\right\}
$$
and
$$
\begin{array}{l}
{\cal B}_E{\cal D}_E=\left\{ B(s)\times D(s): B(s)\times
D(s)\in{\cal B}_E(s)\times {\cal D}_E(s)\right\}
\end{array}.
$$
By  Proposition 1 and Lemma 2, we get the main result of this paper.

{\bf Theorem 1}\ \ ${\cal M}(s)$ is stable for all $B(s)\in{\cal
B}(s),D(s)\in{\cal D}(s)$ if and only if
${\cal M}(s)$ is stable for all $B(s)\times D(s)\in {\cal BD}\cap {\cal B}_E{\cal D}_E$.\\
Proof\ \ By Proposition 1, ${\cal M}(s)$ is stable for all $B(s)\in{\cal
B}(s),D(s)\in{\cal D}(s)$ if and only if
${\cal M}(s)$ is stable for all $B(s)\times D(s)\in  {\cal B}_E{\cal D}_E$.
For any $n(s)\in {\cal M}(s)$ with $B(s)\times D(s)\in  {\cal B}_E{\cal D}_E$,
fixing the rows from the second to the $n$-th and applying Lemma 2 to the first row, we have
$$
\mbox{${\cal M}(s)$ is stable for all $B(s)\in{\cal
B}(s),D(s)\in{\cal D}(s)$ if and only if}
$$
$$
\mbox{${\cal M}(s)$ is stable for all $B(s)\times D(s)\in {\cal BD}_1\cap {\cal B}_E{\cal D}_E$}
$$
where ${\cal BD}_1\subset {\cal B}(s)\times {\cal D}(s)$ satisfies (\ref{4}) when $i=1$.
Repeating the same procedure for the remaining rows consecutively, we get the conclusions.

{\bf Remark 1}\ \ Theorem 1 shows that, the criterion for the robust stability of interval
polynomial matrices can be greatly simplified. The dimension of the minimal testing set is reduced to
$n$.

{\bf Remark 2}\ \ Theorem 1  establishes the minimal dimensional testing criteria for robust stability
of multivariable control systems. By making use of the uncertainty information of  ${\cal B}(s)$ and ${\cal D}(s)$,
i.e., not only making use of the uncertainty information in ${\cal B}(s)$ and ${\cal D}(s)$ separately (this is exactly
what Kamal and Dehleh did), but also making use of the uncertainty  information
collectively, we  are able to improve the robust stability  results given by Kamal and Dehleh.

{\bf Remark 3}\ \ Theorem 1 allows us to improve the main results
of Kamal and Dehleh.

\section{Applications}
\par\indent

Based on Proposition 1, Kamal et al.(1996) address several robustness
problems as follows.

{\bf Proposition 2}($^{\cite{kama}}$)\ \ Let $P(S)=B(S)D^{-1}(s)$ and $K(s)=A(s)C^{-1}(s)$ be
right coprime factorization and left coprime factorization of $P(s)$ and $K(s)$ respectively over the ring of
polynomial matrices. Let $B(s)\in {\cal B}(s)$ and $D(s)\in{\cal D}(s)$, and let $K(s)$ be  fixed.
Then the closed loop system in Fig. 1 is stable for all $B(s)\in{\cal B}(s),D(s)\in{\cal D}(s)$ if and only if it is stable
for all $B(s)\in {\cal B}_E(s),D(s)\in{\cal D}_E(s)$.\\
\unitlength=0.5mm
\begin{picture}(200,70)(-60,0)
\put(5,50){\vector(1,0){15}} \put(25,50){\circle{10}}
\put(30,50){\line(1,0){35}} \put(65,45){\framebox(15,10){P(s)}}
\put(80,50){\line(1,0){35}} \put(115,50){\vector(1,0){15}}
\put(115,50){\line(0,-1){20}}
\put(65,25){\framebox(15,10){K(s)}}
\put(115,30){\line(-1,0){35}}
\put(65,30){\line(-1,0){40}}
\put(25,30){\vector(0,1){15}} \put(27,40){\line(1,0){5}}
\put(20,10){Fig 1: \ \  Feedback Systems}
\end{picture}
\par\indent

Let $H_\infty^{n\times m}({\cal C}_+)$ be the space of matrix valued functions $F(s)$ that are analytic
in ${\cal C}_+$ and bounded on the $j\omega$-axis with the norm $\|F\|_\infty=\sup_{\omega\in{\cal R}}
\sigma_{\max}(F(j\omega))$.

{\bf Proposition 3}($^{\cite{kama}}$)\ \ Let ${\cal G}(s)\subset H_\infty^{n\times n}({\cal C}_+)$ be a family of proper
interval transfer function matrices. Let $G(s)=B(s)D^{-1}(s)$ be a right coprime description of $G(s)$ over the ring
of polynomial matrices, with $D(s)$ column reduced; then \\
1)\ \
$\|G\|_\infty<1$ for all $B(s)\in{\cal B}(s),D(s)\in{\cal D}(s)$
if and only if $\|G\|_\infty<1$ for all $B(s)\in {\cal B}_E(s),D(s)\in{\cal D}_E(s)$.
2)\ \
$G(s)$ is strictly positive real for all $B(s)\in{\cal B}(s),D(s)\in{\cal D}(s)$
if and only if $G(s)$ is strictly positive real for all $B(s)\in {\cal B}_E(s),D(s)\in{\cal D}_E(s)$.

A memoryless nonlinearity $\Delta_s:[0,\infty)\times {\cal R}^n\rightarrow {\cal R}^n$ is said to satisfy
a sector condition if for all $y\in{\cal R}^n$ and $t\geq 0$,
\begin{equation}
[\Delta_s(t,y)-K_1y]^T[\Delta_s(t,y)-K_2y]\geq 0
\end{equation}
for some real matrices $K_1$ and $K_2$, where $K=K_2-K_1$ is a symmetric positive definition matrix.

{\bf Proposition 4}($^{\cite{kama}}$)\ \ $G(s)=B(s)D^{-1}(s)$ is a family of interval proper stable systems, where
$B(s)\in{\cal B}(s),D(s)\in{\cal D}(s)$, and $\Delta_s$ is a set of memoryless nonlinearities that satisfy
the sector condition globally with $K>0$, then the system is globally stable if there is $\eta\geq 0$ such that
$G(\infty)=0$ and
$$
Z(j\omega)=I+(1+j\omega\eta)KG(j\omega)>0,\ \ \forall \omega\in{\cal R}
$$
for all $B(s)\in {\cal B}_E(s),D(s)\in{\cal D}_E(s)$.
\par\indent

By resort to Theorem 1, all testing sets above can be reduced to
lower dimensional ones.

{\bf Theorem 2}\ \ Under the same assumptions as in Proposition 2,
 the closed loop system in Fig. 1 is stable for all $B(s)\in{\cal B}(s),D(s)\in{\cal D}(s)$ if and only if it is stable
for all $B(s)\times D(s)\in {\cal BD}\cap {\cal B}_E{\cal D}_E$.

{\bf Theorem 3}\ \ Under the same assumptions as in Proposition 3, we have \\
1)\ \
$\|G\|_\infty<1$ for all $B(s)\in{\cal B}(s),D(s)\in{\cal D}(s)$
if and only if $\|G\|_\infty<1$ for all $B(s)\times D(s)\in {\cal BD}\cap {\cal B}_E{\cal D}_E$.
2)\ \
$G(s)$ is strictly positive real for all $B(s)\in{\cal B}(s),D(s)\in{\cal D}(s)$
if and only if $G(s)$ is strictly positive real for all $B(s)\times D(s)\in {\cal BD}\cap {\cal B}_E{\cal D}_E$.

{\bf Theorem 4}\ \ Under the same assumptions as in Proposition 4,  the system is globally stable if there is $\eta\geq 0$ such that
$G(\infty)=0$ and
$$
Z(j\omega)=I+(1+j\omega\eta)KG(j\omega)>0,\ \ \forall \omega\in{\cal R}
$$
for all $B(s)\times D(s)\in {\cal BD}\cap {\cal B}_E{\cal D}_E$.

\section{Example}
\par\indent

Recall the two-link planar manipulator considered by Kamal and Dehleh.
The characteristic polynomial of the closed loop system is
\begin{equation}\label{ex}
T(\epsilon)=\det \left(M(\theta_d)s^3+K_ds^2+K_p s+K_r\right)
\end{equation}
where $\theta_d\in[0,\pi/2]$ is the joint angle and
$$
\begin{array}{l}
M(\theta_d)=\left(\begin{array}{cc}
3+2\cos(\theta_d)&1+\cos(\theta_d)\\
1+\cos(\theta_d)&1
\end{array}\right)\\
K_d=\left(\begin{array}{cc}
k_{d11}&k_{d12}\\
k_{d21}&k_{d22}\end{array}\right)\\
K_p=\left(\begin{array}{cc}
k_{p11}&k_{p12}\\
k_{p21}&k_{p22}
\end{array}\right)\\
K_r=\left(\begin{array}{cc}
k_{r11}&k_{r12}\\
k_{r21}&k_{r22}
\end{array}\right)
\end{array}
$$
\begin{equation}\label{in}
\begin{array}{c}
k_{d11}\in[6.07-6.07\epsilon,6.07+6.07\epsilon]\\
k_{d12}\in[2.22-2.22\epsilon,2.22+2.22\epsilon]\\
k_{d21}\in[2.22-2.22\epsilon,2.22+2.22\epsilon]\\
k_{d22}\in[1.62-1.62\epsilon,1.62+1.62\epsilon]\\
k_{p11}\in[6.12-6.12\epsilon,6.12+6.12\epsilon]\\
k_{p12}\in[2.24-2.24\epsilon,2.24+2.24\epsilon]\\
k_{p21}\in[2.24-2.24\epsilon,2.24+2.24\epsilon]\\
k_{p22}\in[1.64-1.64\epsilon,1.64+1.64\epsilon]\\
k_{r11}\in[5.11-5.11\epsilon,5.11+5.11\epsilon]\\
k_{r12}\in[1.87-1.87\epsilon,1.87+1.87\epsilon]\\
k_{r21}\in[1.87-1.87\epsilon,1.87+1.87\epsilon]\\
k_{r22}\in[1.37-1.37\epsilon,1.37+1.37\epsilon].
\end{array}
\end{equation}
Since for any fixed $\theta_d\in[0,\pi/2]$, the inertia matrix
$$
M(\theta_d)=
\left(\begin{array}{cc}
1&(1+\cos(\theta_d))\\
(-1+\cos(\theta_d))&1
\end{array}\right)
\left(\begin{array}{cc}
1&0\\
2&1
\end{array}\right)
$$
In what follows, we take the uncertain inertia matrix as
$$
M(s)=B(s)A(s)
$$
where $A(s)=\left(\begin{array}{ll}1&0\\
2&1
\end{array}\right)$, $B(s)\in{\cal B}(s)=\{(b_{ij}(s))_{2\times 2}:b_{ij}(s)\in B(ij)\}$
and
$$
\begin{array}{ll}
B(11)=s^3 & B(12)=b_{12}s^3,\; b_{12}\in[1,2]\\
B(21)=b_{21}s^3,\; b_{21}\in[-1,0]\ \ & B(22)=s^3
\end{array}
$$
Let ${\cal D}(s)=K_d s^3+K_ps^2+K_r$, then
$$
{\cal D}(s)=\left\{(d_{ij}(s)):d_{ij}(s)\in D(ij)\right\},
$$
where $D(ij)=\left\{d_{ij}(s)=k_{dij}s^2+k_{pij}s+k_{rij}:k_{dij},k_{pij} \mbox{ and } k_{rij}
\mbox{ are given in (\ref{in})}\right\}$.
Thus, the characteristic polynomial of the uncertain system is
\begin{equation}
T(\epsilon)=\det\left(B(s)A(s)+D(s)\right),\ \ B(s)\in{\cal B}(s),D(s)\in{\cal D}(s)
\end{equation}
where $A(s)$ is a fixed matrix and $B(s),D(s)$ are interval polynomial matrices.
\par\indent

The Kharitonov polynomial vertex sets associated with $B(ij)$ and $D(ij)$ are as follows.
$$
\begin{array}{c}
\begin{array}{cccc}
{\cal B}_V(11)=\{s^3\}=B(11), & {\cal B}_V(12)=\{s^3,2s^3\}, & {\cal B}_V(21)=\{-s^3,0\}, & {\cal B}_V(22)=\{s^3\}=B(22),
\end{array}\\
{\cal D}_V(ij)=\left\{\begin{array}{cc}
k_{ij}^{(1)}=K_{dij}^U s^2+K_{pij}^L s+K_{rij}^L,& k_{ij}^{(2)}=K_{dij}^U s^2+K_{pij}^U s+K_{rij}^L,\\
k_{ij}^{(3)}=K_{dij}^L s^2+K_{pij}^L s+K_{rij}^U,& k_{ij}^{(4)}=K_{dij}^L s^2+K_{pij}^U s+K_{rij}^U,
\end{array}\right\}
\end{array}
$$
where the upper bound and the lower bound of $K_{dij}$ are denoted as $K_{dij}^U,K_{dij}^L$.
And the Kharitonov edge sets associated with $B(ij)$ and $D(ij)$ are  as follows.
$$
\begin{array}{c}
\begin{array}{cc}
{\cal B}_E(11)=\{s^3\}=B(11), & {\cal B}_E(12)=\{\lambda s^3+(1-\lambda)\cdot 2s^3:\lambda\in[0,1]\},\\
{\cal B}_E(21)=\{-\lambda s^3:\lambda\in[0,1]\}, & {\cal B}_E(22)=\{s^3\}=(22)
\end{array}\\
{\cal D}_E(ij)=\left\{
\begin{array}{c}
\lambda k_{ij}^{(1)}+(1-\lambda)k_{ij}^{(2)},\\
\lambda k_{ij}^{(2)}+(1-\lambda)k_{ij}^{(4)},\\
\lambda k_{ij}^{(4)}+(1-\lambda)k_{ij}^{(3)},\\
\lambda k_{ij}^{(3)}+(1-\lambda)k_{ij}^{(1)}.
\end{array}\right\}
\end{array}
$$
By Proposition 1 given by Kamal and Dehleh, for any $\epsilon$, $T(\epsilon)$ is robustly stable
if and only if $T(\epsilon)$ is stable for all $B(s)\times D(s)\in {\cal B}_E(s)\times {\cal D}_E(s)$, where
$$
\begin{array}{c}
{\cal B}_E(s)\times {\cal D}_E(s)=U_1\cup U_2\\
U_1= \left(\begin{array}{cc}
B(11)&B_E(12)\\
B_E(21)&B(22)
\end{array}\right)\times
\left(\begin{array}{cc}
D_V(11)&D_E(12)\\
D_E(21)&D_V(22)
\end{array}\right)\\
U_2= \left(\begin{array}{cc}
B(11)&B_E(12)\\
B_E(21)&B(22)
\end{array}\right)\times
\left(\begin{array}{cc}
D_E(11)&D_V(12)\\
D_V(21)&D_E(22)
\end{array}\right)
\end{array}.
$$
So, the robust stability of $T(\epsilon)$ is reduced to $2\cdot 4^6$ robust stability problems involving
determinants with four parameters. In the case, the testing set is a four-dimensional set. In the sequel, we apply Theorem 1 to simplify
this problem further. That is to say,
$T(\epsilon)$ is robustly stable
if and only if $T(\epsilon)$ is stable for all
$B(s)\times D(s)\in  ({\cal B}{\cal D})\cap {\cal B}_E{\cal D}_E$, where
$$
\begin{array}{c}
({\cal B}{\cal D})\cap {\cal B}_E{\cal D}_E=V_1\cup V_2\cup V_3\cup V_4\cup V_5\cup V_6\cup V_7\\
V_1=\left(\begin{array}{cc}
B(11)&B_E(12)\\
B_E(21)&B(22)
\end{array}\right)\times
\left(\begin{array}{cc}
D_V(11)&D_V(12)\\
D_V(21)&D_V(22)
\end{array}\right)\\
V_2=\left(\begin{array}{cc}
B(11)&B_E(12)\\
B_V(21)&B(22)
\end{array}\right)\times
\left(\begin{array}{cc}
D_E(11)&D_V(12)\\
D_V(21)&D_V(22)
\end{array}\right)\\
V_3=\left(\begin{array}{cc}
B(11)&B_E(12)\\
B_V(21)&B(22)
\end{array}\right)\times
\left(\begin{array}{cc}
D_V(11)&D_V(12)\\
D_E(21)&D_V(22)
\end{array}\right)\\
V_4=\left(\begin{array}{cc}
B(11)&B_V(12)\\
B_E(21)&B(22)
\end{array}\right)\times
\left(\begin{array}{cc}
D_V(11)&D_E(12)\\
D_V(21)&D_V(22)
\end{array}\right)\\
V_5=\left(\begin{array}{cc}
B(11)&B_V(12)\\
B_E(21)&B(22)
\end{array}\right)\times
\left(\begin{array}{cc}
D_V(11)&D_V(12)\\
D_V(21)&D_E(22)
\end{array}\right)\\
V_6=\left(\begin{array}{cc}
B(11)&B_V(12)\\
B_V(21)&B(22)
\end{array}\right)\times
\left(\begin{array}{cc}
D_E(11)&D_V(12)\\
D_V(21)&D_E(22)
\end{array}\right)\\
V_7=\left(\begin{array}{cc}
B(11)&B_V(12)\\
B_V(21)&B(22)
\end{array}\right)\times
\left(\begin{array}{cc}
D_V(11)&D_E(12)\\
D_E(21)&D_V(22)
\end{array}\right)
\end{array}
$$
Clearly, the robust stability of $T(\epsilon)$ is reduced to $7\cdot 4^6$ robust stability problems involving
determinants with two parameters. In this case, the minimal testing set is two-dimensional.

{\bf Remark 4:}\ \ A test with a parameter could involve infinite tests without parameter, such as an edge
means infinite points. The significance of our result is that it is the minimal dimension test for such problem as
(\ref{1}). Furthermore, we consider the uncertainty in inertia matrix
$M(\theta_d)$ in Example, which did not take into consideration in Ref. \cite{kama}.

\section{Conclusions}
\par\indent

In this paper, interval transfer function matrices are studied.
For a family of multivariable control systems having interval transfer functions  matrices,
a lower dimensional robust stability criterion has been established. This result improves the
previous criterion obtained by Kamal and Dahleh, and can be used to establish improved lower
dimensional robust performance criteria as well.


\begin{thebibliography}{19}
{\small
\bibitem{khar}V. L. Kharitonov. Asymptotic stability of an equilibrium
  position of a family of systems of linear differential equations,
Differential'nye Uravneniya, vol.14, 2086-2088, 1978.

\bibitem{bart}A. C. Bartlett, C. V. Hollot and L. Huang. Root locations of
an entire polytope of polynomials: It suffices to check the edges,
Mathematics of Control, Signals, and Systems, vol.1, 61-71, 1988.

\bibitem{kama}F. Kamal and M. Dahleh. Robust stability of multivariable interval control systems,
Int. J. Control, vol.64, 807-828, 1996.

\bibitem{wang}L. Wang and L. Huang. Vertex results for uncertain
systems, Int. J. Systems Science, vol.25, 541-549, 1994.

\bibitem{long}L. Wang,  Z. Z. Wang and W. S. Yu. Stability of polytopic polynomial
matrices. {\it Proc. American Conrol Conference}, Arlington,
Virginia, June: 4695-4696,2001.

\bibitem{chap1}H. Chapellat and S. P. Bhattacharyya. A generalization of Kharitonov's
theorem: robust stability of interval plants. IEEE Trans. on
Automatic Control, vol.34, 306-311, 1989.

\bibitem{bias}S. Bialas. A necessary and sufficient condition for the
stability of interval matrices, Int. J. Control, vol.37, 717-722, 1983.

\bibitem{barm}B. R. Barmish and C. V. Hollot. Counter-example to a
recent result on the stability of interval matrices by S. Bialas,
Int. J. Control, vol.39, 1103-1104, 1984.

\bibitem{koka}H. Kokame and T. Mori. A Kharitonov-like theorem for
interval polynomial matrices, Systems and Control Letters,
vol.16, 107-116, 1991.

\bibitem{acke}{J. Ackermann,
{\it Robust Control: Systems with Uncertain Physical Parameters.}
Springer-Verlag, London, 1993.}

\bibitem{barm1}B. R. Barmish, New tools for robustness analysis,
Proc. of IEEE Conf. on Decision and Control, 1-6, 1988.
}
\end{thebibliography}
\end{document}